\documentclass[12pt,a4paper]{article}
\usepackage{a4wide}
\usepackage{amsthm} 
\usepackage{amsmath}
\usepackage{amssymb} 
\usepackage{latexsym} 
\usepackage{array} 
\usepackage[T1]{fontenc} %
\usepackage{fancyhdr} 
\usepackage{xspace} %
\usepackage[dvips]{epsfig} 
\usepackage{graphicx} 
\usepackage{enumerate}
\renewcommand{\theequation}{\arabic{section}.\arabic{equation}}

\newcommand{\R}{{\mathbb R}}
\newcommand{\N}{{\mathbb N}}
\newcommand{\on}{\operatorname}
\def\Om{\Omega}
\def\om{\omega}
\def\p{\partial}
\def\ds{\displaystyle}
\def\lf{\left}
\def\rg{\right}
\def\ep{\varepsilon}
\def\la{\lambda}
\def\ti{\tilde}
\newtheorem{definition}{Definition}[section]
\newtheorem{bem}[definition]{Remark}

\newtheorem{theorem}[definition]{Theorem}
\newtheorem{lemma}[definition]{Lemma}

\newcommand{\clearemptydoublepage}%

  {\newpage{\pagestyle{empty}\cleardoublepage}}
\author{\bf{Tobias Lamm}\footnote{Departement Mathematik, ETH Z\"urich, 8092 Z\"urich, Switzerland, tobias.lamm@math.ethz.ch} and \bf{Tristan Rivi{\`e}re}\footnote{Departement Mathematik, ETH Z\"urich, 8092 Z\"urich, Switzerland, tristan.riviere@math.ethz.ch}}

\title{\bf{Conservation laws for fourth order systems in four dimensions}\rm}
\begin{document}
\maketitle
\begin{abstract}
Following an approach of the second author \cite{riviere06} for conformally invariant variational problems in two dimensions, we show in four dimensions the existence of a conservation law for fourth order systems, which includes both intrinsic and extrinsic biharmonic maps. With the help of this conservation law we prove the continuity of weak solutions of this system. Moreover we use the conservation law to derive the existence of a unique global weak solution of the extrinsic biharmonic map flow in the energy space. 
\end{abstract}
\section{Introduction}
In \cite{riviere06} the second author found a procedure to rewrite second order elliptic systems of the form
\begin{equation}
\label{z1}
-\Delta u=\Omega\cdot \nabla u
\end{equation}
in divergence form. Here $u$ is a map into ${\R}^m$ and $\Omega=\Om_k\, \p_{x_k}$ is a  vectorfield tensored with $m\times m-$anti-symmetric matrices.
The procedure consists of finding a map $A$ which takes values in the space of $m\times m-$invertible matrices and a two-vectorfield tensored with $m\times m-$matrices $B=B_{kl}\ \p_{x_k}\wedge \p_{x_l}$ satisfying
\begin{equation}
\label{z2}
\nabla A -A\ \Omega+curl\ B=0
\end{equation}
where $curl B$ is the matrix valued vector field $(\sum_k\p_{x_k} B_{kl})\ \p_{x_l}$. Once the existence of ``regular enough'' $A$ and $B$ satisfying (\ref{z2}) has been established, one observes that the system (\ref{z1}) is equivalent to the following {\it conservation law}
\begin{equation}
\label{z3}
div \lf( A\nabla u+B\cdot\nabla u\rg)=0
\end{equation}
where the latter equation means in coordinates $\p_{x_k}(A\p_{x_k}u+B_{kl}\p_{x_l}u)=0$ with implicit summation convention.
It was shown in \cite{riviere06} that the above procedure works successfully in two dimensions for solutions $u \in W^{1,2}$ and for $\Omega \in L^2$. Writing the equation in the form (\ref{z3}) permits to prove, without many additional efforts, results such as the continuity of $W^{1,2}$ solutions to (\ref{z1}) in two dimensions for a connection field $\Om \in L^2$ or the sequentially weak compactness in $W^{1,2}$ of these solutions, etc. It was also observed in \cite{riviere06} that every Euler Lagrange equation of an elliptic conformally invariant Lagrangian with quadratic growth in two dimensions can be written in the form (\ref{z1}). In \cite{rivierestruwe06} Struwe and the second author used the existence of an "almost" conservation law, which was motivated by the conservation law \eqref{z3}, to give a new proof for the partial regularity of harmonic maps, and generalizations thereof, in higher dimensions.
\newline
Working now in four dimensions and replacing the Laplacian by the Bilaplacian it is then natural to ask how far the previous results can be extended to this new setting. We shall use the following notation : $L^2\cdot W^{1,2}(B^4,M(m))$ denotes the space of linear combination of products of an $L^2$ map with an $W^{1,2}$ map from $B^4$ into the space of $m\times m -$matrices $M(m)$.
Our main result in this work is the following   
\begin{theorem}
\label{th-1}
Let $u \in W^{2,2}(B^4,{\R}^m)$ satisfy the equation
\begin{equation}
\label{z4}
\Delta^2 u=\Delta\lf(V\cdot\nabla u\rg)+ div\lf(w\ \nabla u\rg)+W\cdot \nabla u,
\end{equation}
where the ``potentials'' $V$ and $w$ are in $W^{1,2}(B^4,M(m)\otimes\wedge^1{\R}^4)$ respectively $L^2(B^4,M(m))$ and the potential $W \in W^{-1,2}(B^4,M(m)\otimes\wedge^1{\R}^4)$ can be decomposed in the following way: there exists $\om \in L^2(B^4,so(m))$ and $F \in L^2\cdot W^{1,2}(B^4,M(m)\otimes\wedge^1{\R}^4 )$ such that
\begin{equation}
\label{z5}
W=\nabla\om+F.
\end{equation}
Then $u$ is continuous in $B^4$.
\end{theorem}
\begin{bem}
The Theorem remains true if we only assume that $F \in L^{\frac{4}{3},1}(B^4,M(m)\otimes\wedge^1{\R}^4 )$ (for a definition of the Lorentz space $L^{\frac{4}{3},1}$ see section $2.2.1$).
\end{bem}
As in \cite{riviere06} for solutions to (\ref{z1}), the previous result is based on the discovery of a conservation law satisfied by solutions of (\ref{z4}). Precisely we establish in the present work the following theorem.
\begin{theorem}
\label{th-2}
Let $V$ and $w$ be in $W^{1,2}(B^4,M(m)\otimes\wedge^1{\R}^4)$ respectively $L^2(B^4,M(m))$ and let $W \in W^{-1,2}(B^4,M(m)\otimes\wedge^1{\R}^4)$ satisfy (\ref{z5}) for $F \in L^2\cdot W^{1,2}(B^4,M(m)\otimes\wedge^1{\R}^4 )$. Let $A \in L^\infty\cap W^{2,2}(B^4,Gl(m))$ and $B \in L^2(B^4, M(m)\otimes \wedge^2 \R^4)$ satisfy the linear equation
\begin{equation}
\label{z6}
\nabla\Delta A+\Delta A\ V-\nabla A\ w+A\ W=curl\, B,
\end{equation}
where $curl\, B:=\sum_l\p_{x_l} B_{lk}\ \p_{x_k}$. Then $u$ solves the equation (\ref{z1}) if and only if it satisfies the conservation law
\begin{equation}
\label{z7}
\begin{array}{l}
div[\ds\nabla(A\Delta u)-2\nabla A\ \Delta u+\Delta A\ \nabla u-A\, w\ \nabla u\\[5mm]
 \quad\quad\ds+\nabla A\lf(V\cdot\nabla u\rg)-A\ \nabla(V\cdot\nabla u)
-B\cdot\nabla u]=0,
\end{array}
\end{equation}
where we use the notation $B\cdot\nabla u:=B_{kl}\p_{x_l}u\ \p_{x_k}$.
\end{theorem} 
\begin{bem}\label{rem}
We expect similar Theorems to remain true for general even order elliptic systems of the type \eqref{z4}. 
\end{bem}
Theorem~\ref{th-1} will be a consequence of the previous result and the following existence result for $A$ and $B$.
\begin{theorem}
\label{th-3}
There exists $\ep(m)>0$ such that the following holds: Let $V$ and $w$ respectively be in $W^{1,2}(B^4_2,M(m)\otimes\wedge^1{\R}^4)$ and in $L^2(B^4_2,M(m))$ and
let $W=\nabla \om+F \in W^{-1,2}(B^4_2,M(m)\otimes\wedge^1{\R}^4)$ satisfy (\ref{z5}) for $F \in L^2\cdot W^{1,2}(B^4_2,M(m)\otimes\wedge^1{\R}^4 )$. Assume that
\begin{equation}
\label{z8}
\|V\|_{W^{1,2}}+\|w\|_{L^2}+\|\om\|_{L^2}+\|F\|_{L^2\cdot W^{1,2}}<\ep(m),
\end{equation}
then there exist $A \in L^\infty\cap W^{2,2}(B^4,Gl(m))$ and  $B \in W^{1,4/3}(B^4,M(m) \otimes \wedge^2{\R}^4)$ satisfying (\ref{z6}) on $B^4$ with the following estimate
\begin{equation}
\label{z8a}
\begin{array}{l}
\ds\|A\|_{W^{2,2}}+\|dist(A,SO(m))\|_{L^\infty}+\|B\|_{W^{1,{4}/{3}}}\\[5mm]
\ds\quad\quad\le C\ (\lf\|V\|_{W^{1,2}}+\|w\|_{L^2}+\|\om\|_{L^2}+\|F\|_{L^2\cdot W^{1,2}}\rg).
\end{array}
\end{equation}
\end{theorem}
The previous theorems apply to classical critical fourth order non-linear elliptic systems in four dimensions such as intrinsic and the extrinsic biharmonic maps into Riemannian manifolds. 
\begin{theorem}
\label{th-4}
Let $N^n$ be a closed submanifold of ${\R}^m$. Let $u$ in $W^{2,2}(B^4,N^n)$ be a critical point in $W^{2,2}(B^4,N^n)$ of the functional
\begin{equation}
\label{z9}
{\mathcal B}_{ext}(u)=\int_{B^4}|\Delta u|^2\ dx,
\end{equation}
or of the functional
\begin{equation}
\label{z10}
{\mathcal B}_{int}(u)=\int_{B^4}|(\Delta u)^T|^2\ dx,
\end{equation}
where $(\Delta u)^T$ is the tension field of $u$ : the projection of $\Delta u$ on the tangent space $T_{u(x)}N$ to $N$ at $u(x)$. Then $u$ satisfies an equation of the form (\ref{z1}) where $V$, $w$ and $W$ satisfy the assumptions of Theorem~\ref{th-1}.
\end{theorem}
\begin{bem}
Another example of equations which can be written in the form \eqref{z4} is given by
\begin{equation}
\label{z10a}
-\Delta^2 u= Q(x,u,Du),
\end{equation}
where $Q$ satisfies 
\begin{equation}
\label{z10b}
|Q(x,z,p)|\le c |p|^4 \ \ \ \forall\ \ \ (x,z,p)\in \Omega \times \R^n \times \R^{4n}.
\end{equation} 
In contrast to the corresponding two dimensional problem (see \cite{frehse73}) this implies in particular that a quartic nonlinearity in the gradient is not critical for fourth order systems in four dimensions. The regularity of solutions of \eqref{z10a} was first studied by Wang \cite{wang1}.
\end{bem}
Therefore, Theorem~\ref{th-2} and Theorem~\ref{th-3} permit to write the extrinsic and intrinsic biharmonic map equation into arbitrary target manifolds in divergence form and this yields a new and simple proof of the regularity of these maps in four dimensions (for previous results see \cite{chang99}, \cite{strzelecki03}, \cite{wang1} and \cite{wang3}). The illustration of the above result by giving the explicit values of $V$, $w$, $W$, $A$ and $B$ in the particular case where $N^n$ is the unit sphere $S^n$ of ${\R}^{n+1}$ is instructive: a classical computation (see \cite{chang99}, \cite{strzelecki03} and \cite{wang3}) tells us that $u$ is an extrinsic biharmonic map (i.e. a $W^{2,2}$ critical point of (\ref{z9}) for perturbations in $W^{2,2}(B^4,S^n)$) if and only if it satisfies (\ref{z1}) for $V$, $w$ and $W$ given by
\begin{equation}
\label{z11}
\lf\{
\begin{array}{l}
V^{ij}=u^i\nabla u^j-u^j\nabla u^i\\[5mm]
w^{ij}=div\ V^{ij}-2|\nabla u|^2\delta^{ij}\\[5mm]
\ds W^{ij}=\nabla(div\ V^{ij})+2\lf[\Delta u^i\ \nabla u^j-\Delta u^j\ \nabla u^i\rg]
\end{array}
\rg.
\end{equation}
Observe that $div\ W=0$. This implies that there exists $B \in L^2(B^4,M(m)\otimes\wedge^2{\R}^4)$ such that $W=curl\, B$. Then $u$ satisfies (\ref{z7}) with
\begin{equation}
\label{z12}
\lf\{
\begin{array}{l}
A=I_m\\[5mm]
B\quad\mbox{ s.t.}\quad W=curl\, B.
\end{array}
\rg.
\end{equation}
This conservation law was known in the particular case where the target $N^n$ is the standard round sphere $S^n \hookrightarrow \R^n$. The main contribution of the present work is to show that this conservation law is stable and keeps existing while changing the target to an arbitrary one. Similarly $u$ is an intrinsic biharmonic map (i.e. a $W^{2,2}$ critical point of (\ref{z10}) for perturbations in $W^{2,2}(B^4,S^n)$) if and only if it satisfies (\ref{z1}) for $V$, $w$ and $W$ given by
\begin{equation}
\label{z13}
\lf\{
\begin{array}{l}
V^{ij}=u^i\nabla u^j-u^j\nabla u^i\\[5mm]
w^{ij}=div\ V^{ij}\\[5mm]
\ds W^{ij}=\nabla(div\ V^{ij})+2\lf[\Delta u^i\ \nabla u^j-\Delta u^j\ \nabla u^i+|\nabla u|^2\lf(u^i\ \nabla u^j-u^j\ \nabla u^i\rg)\rg]
\end{array}
\rg.
\end{equation}
Observe that in this case also $div\ W=0$. Again there exists $B \in L^2(B^4,M(m)\otimes\wedge^2{\R}^4)$ such that $W=curl\, B$. Then $u$ satisfies (\ref{z7}) with
\begin{equation}
\label{z14}
\lf\{
\begin{array}{l}
A=I_m\\[5mm]
B\quad\mbox{ s.t.}\quad W=curl\, B
\end{array}
\rg.
\end{equation} 
In the second part of the paper we use the conservation law \eqref{z7} to show the global existence of a unique weak solution of the extrinsic biharmonic map flow (i.e. the gradient flow for the energy $\mathcal{B}_{ext}$) in four dimensions in the energy space. More precisely we consider a smooth, compact, four-dimensional Riemannian manifold $M$ without boundary and we study the parabolic system,
\begin{equation}
\label{flow}
\begin{array}{rl}
\ds \partial_t u =& -\Delta^2 u -\Delta (V\cdot \nabla u)-\text{div} (w\nabla u)-W\cdot \nabla u \ \ \ \text{in}\ \ \ \mathcal{D}'(M \times (0,T)),   \\[5mm]
\ds u(\cdot,0)=& u_0,
\end{array}
\end{equation} 
where $V$, $w$ and $W$ are as in \eqref{extbih} and $u_0 \in W^{2,2}(M,N)$. As mentioned above we prove the following 
\begin{theorem}\label{flow1}
For all $u_0 \in W^{2,2}(M,N)$ there exists a unique global weak solution $u\in H^1([0,T],L^2(M,N)) \cap L^2([0,T],W^{2,2}(M,N))$ of \eqref{flow} with non-increasing energy $\mathcal{B}_\text{ext}$. 
\end{theorem} 
In \cite{lamm03b} the first author proved the longtime existence of a smooth solution of \eqref{flow} assuming a smallness condition on the initial energy. Gastel \cite{gastel04} and Wang \cite{wang06} proved the existence of a unique global weak solution (which is smooth away from finitely times) of \eqref{flow} and higher order generalisations of this flow. 
\newline
The corresponding theorem for the harmonic map flow was proved by Freire \cite{freire95}, \cite{freire95b}, \cite{freire95c} following previous work of the second author \cite{riviere93}. Global existence of solutions of gradient flow for $\mathcal{B}_{int}$ has been shown by the first author \cite{lamm05} in the case of a target manifold with non-positive sectional curvature.
\section{Proof of Theorems~\ref{th-1}.--\ref{th-4}.} 
\setcounter{equation}{0}
\subsection{Proof of Theorem~\ref{th-2}.}
Theorem~\ref{th-2} is proved by a direct computation. More precisely we have for general $A$, $B$, $u$, $V$, $w$, $W$ in the spaces given in the statement of the Theorem that
\begin{equation}
\label{z20}
\begin{array}{l}
\ds div\lf[\nabla( A\Delta u)-2\nabla A\ \Delta u+\Delta A\ \nabla u -Aw\ \nabla u+\nabla A\ (V\cdot\nabla u)-A\ \nabla(V\cdot\nabla u)\rg]\\[5mm]
\quad=A\ \lf[\Delta^2 u-\Delta\lf(V\cdot\nabla u\rg)- div\lf(w\ \nabla u\rg)-W\cdot \nabla u\rg]\\[5mm]
\quad +\lf[\Delta \nabla A+\Delta A\ V-\nabla A\ w+A\  W\rg]\cdot\nabla u
\end{array}
\end{equation}
Combining this equation with the fact that 
\[
 div(B\cdot\nabla u)=\p_{x_k}( B_{kl}\p_{x_l}u)=\p_{x_k}B_{kl}\ \p_{x_l}u=curl\, B\cdot\nabla u
\]
(we are using the fact that $B$ takes values in $M(m) \otimes \wedge^2{\R}^4$  and thus : $B_{kl}=-B_{lk}$) we deduce that (\ref{z7}) holds if and only if $u$ is a solution of equation (\ref{z4}). This proves Theorem~\ref{th-2}. 
\subsection{Proof of Theorem~\ref{th-3}.}

\subsubsection{Some results on Lorentz spaces.}
The Lorentz spaces will play an important role in the proof of Theorem~\ref{th-3} and we recall some classical facts about these spaces which where proved in \cite{oneil63}, \cite{peetre63}, \cite{poornima83} and also exposed in \cite{helein02}, \cite{hunt66}, \cite{stein71} and \cite{tartar98}. Let $f$ be a measurable function on $\Om$ a domain in ${\R}^k$ and denote by $f^\ast(t)$ the equimeasurable decreasing rearrangement of $f$ which is a function on ${\R}_+ $ satisfying
\[
\lf|x\in \Om\ ;\ |f|(x)>\la\rg|=\lf|t\in{\R}_+\ ;\ f^\ast(t)>\la\rg|
\]
For $f$ measurable in $\Om$ we introduce
\[
f^{\ast\ast}(t)=\frac{1}{t}\int_0^tf^\ast(s)\ ds
\]
and, for every $1\le p\le \infty$ and $1\le q\le\infty$, we define the following functional which happens to be a norm on the set of measurable functions for which it is finite
\[
\|f\|_{L^{p,q}}=
\lf\{
\begin{array}{l}
\ds\lf(\int_0^\infty(t^{1/p}f^{\ast\ast}(t))^q\ \frac{dt}{t}\rg)^{1/q}\quad\quad\mbox{ if }1\le q<\infty\\[5mm]
\ds\sup_{t>0}t^{1/p}f^{\ast\ast}(t)\quad\quad\mbox{ if }q=\infty.
\end{array}
\rg.
\]
The set of functions $f$ for which $\|f\|_{L^{p,q}}$ is finite is a Banach space for this norm and it is called Lorentz $L^{p,q}$ space. It happens that $L^{p,p}$ is the standard $L^p$ space, that, if $\Om$ is bounded, $L^{p,q}(\Om)$ embeds into $L^{p',q'}(\Om)$ if and only if $p>p'$ or $q'\ge q$ in the case that $p=p'$. We shall use the following classical results on Lorentz spaces: first of all $L^{p,q}\cdot L^{p',q'}$ embeds continuously into $L^{r,s}$ if $\frac{1}{p}+\frac{1}{p'}\le 1$ for
\[
\frac{1}{r}=\frac{1}{p}+\frac{1}{p'}\quad\mbox{ and }\quad\frac{1}{s}=\frac{1}{q}+\frac{1}{q'},
\]
and
\[
\|fg\|_{L^{r,s}}\le C\  \|f\|_{L^{p,q}}\ \|g\|_{L^{p',q'}}.
\]
Moreover we will use the fact that Calderon-Zygmund operator maps continuously $L^{p,q}$ into $L^{p,q}$ for $1<p<+\infty$ and that the dual space of $L^{2,1}$ is $L^{2,\infty}$. Finally, the Lorentz-Sobolev-space $W^{m,p,q}(\Om)$ of functions whose first $m$ derivatives are in $L^{p,q}$ embeds into $L^{p^\ast,q}$ where $1/p^\ast=1/p-m/k$.  We shall use the previous results in the following situations: $W^{1,2}(B^4)$ embeds into $L^{4,2}(B^4)$, the product of two functions in $L^{4,2}$ is in $L^{2,1}$, the product of a function in $L^2$ with a function in $W^{1,2}$ is in $L^{4/3,1}$, a function whose Laplacian is in $L^{2,1}(B^4)$ is continuous (see also the proof of Theorem \ref{th-1}), etc.
\subsubsection{Proof of Theorem~\ref{th-3}.}
 Let $V$ and $w$ be in $W^{1,2}(B^4_2,M(m)\otimes\wedge^1{\R}^4)$ respectively $L^2(B^4_2,M(m))$ and let $W=\nabla \om+F \in W^{-1,2}(B^4_2,M(m)\otimes\wedge^1{\R}^4)$ satisfy (\ref{z5}) for $F \in L^2\cdot W^{1,2}(B^4_2,M(m)\otimes\wedge^1{\R}^4 )$. Assume that
\begin{equation}
\label{z8b}
\Delta=\|V\|_{W^{1,2}}+\|w\|_{L^2}+\|\om\|_{L^2}+\|F\|_{L^2\cdot W^{1,2}}<\ep(m),
\end{equation}
where $\ep(m)$ will be chosen small enough later. Using standard elliptic theory we get the existence of $\Om \in W^{1,2}(B^4_2, so(m) \otimes \wedge^1{\R}^4)$ satisfying 
\begin{equation}
\label{z21}
\lf\{
\begin{array}{l}
div(\Om)=-\om\quad\mbox{ in }B^4_2 \\[5mm]
\|\Om\|_{W^{1,2}}\le c\ \|\om\|_{L^2}\le c\ \ep(m)
\end{array}
\rg.
\end{equation}
Then for $\ep(m)$ small enough we can apply Theorem~\ref{Uhlenbeckgauge} in order to get a Coulomb gauge: a 2-vectorfield $\xi=\xi_{kl}\ \p_{x_k}\wedge\p_{x_l}\in W^{2,2}(B^4_2,so(m)\otimes \wedge^2 \R^4)$ (i.e. $\xi_{kl}=(\xi^{ij}_{kl})_{ij}\in so(m)$ and we have $\xi^{ij}_{kl}=-\xi^{ji}_{kl}$ and $\xi^{ij}_{kl}=-\xi_{lk}^{ij})$ and a map $U \in W^{2,2}(B_2^4,so(m))$ such that
\begin{equation}
\label{z22}
\lf\{
\begin{array}{l}
\ds \Om=P\nabla P^{-1}+P\  curl\, \xi\ P^{-1}\\[5mm]
\ds P:=\exp(U)\\[5mm]
\ds  \|U\|_{W^{2,2}(B_2^4)}+\|\xi\|_{W^{2,2}(B_2^4)}\le c\ \|\om\|_{L^2(B_2^4)}.
\end{array}
\rg.
\end{equation}
We calculate
\begin{equation}
\label{z23}
\begin{array}{rl}
\ds\nabla \om&\ds=-\nabla div \Om\\[5mm]
 &\ds=-\nabla\lf(\nabla P\cdot\nabla P^{-1}+P\ \Delta P^{-1}+div(P\ curl\,\xi\ P^{-1})\rg)\\[5mm]
 &\ds= -P\ \nabla\Delta\ P^{-1}+K_1.
 \end{array}
 \end{equation}
Using the fact that 
\[
div(P\ curl\,\xi\ P^{-1})=\nabla P\cdot curl\,\xi\ P^{-1}+P\ curl\,\xi\cdot\nabla P^{-1},
\]
(\ref{z22}) and the properties of Lorentz spaces we mentioned in the previous subsection, we have
\begin{equation}
\label{z24}
\|K_1\|_{L^{4/3,1}}\le c\ \|\om\|^2_{L^2}. 
\end{equation}
We now proceed with some a-priori computations (before having the existence of $A$ and $B$). Consider an element $A \in W^{2,2}\cap L^\infty(B_2^4,Gl(m))$ and an element $P$ as in \eqref{z22}. Denoting $\ti{A}:=A\, P$ we then have
\begin{equation}
\label{z25}
\begin{array}{l}
\ds\nabla\Delta\ti{A}-\nabla\Delta(\ti{A}\,P^{-1})\ P\\[5mm]
\ds= -\ti{A}\ \nabla\Delta P^{-1}\ P-\nabla\ti{A}\ \Delta P^{-1}\ P-\Delta\ti{A}\ \nabla P^{-1}\ P
-2\nabla(\nabla\ti{A}\cdot\nabla P^{-1})\ P.
\end{array}
\end{equation}
Combining this with \eqref{z6} and \eqref{z23} we get
\begin{equation}
\label{z251}
\begin{array}{rl}
\ds \nabla \Delta \ti{A}\ =&-\Delta (\tilde{A}P^{-1})V P+\nabla (\tilde{A}P^{-1})wP-\tilde{A}(-\nabla \Delta P^{-1}+P^{-1}K_1)P+(curl B) P\\[5mm]
\ds &-\ti{A}\ \nabla\Delta P^{-1}\ P-\nabla\ti{A}\ \Delta P^{-1}\ P-\Delta\ti{A}\ \nabla P^{-1}\ P -2\nabla(\nabla\ti{A}\cdot\nabla P^{-1})\ P\\[5mm]
=& -\Delta \ti{A} K_2-\nabla^2\ti{A} K_3+ \nabla \ti{A} K_4 -\ti{A} K_5+(curl B) P,
\end{array}
\end{equation}
where we have the estimate
\begin{equation}
\label{z252}
\begin{array}{l}
||K_2||_{W^{1,2}}+||K_3||_{W^{1,2}}+||K_4||_{L^2}+||K_5||_{L^{\frac{4}{3},1}} \le c \Delta.
\end{array}
\end{equation}
Instead of aiming to solve \eqref{z251} we use a cut-off function $\phi\in C^\infty_0(B^4_2)$, with $\phi=1$ in $B^4$ and $|\nabla^l \phi|\le c$ for all $l\in \N$, to get new maps $\overline{V}=\phi V$, $\overline{w}=\phi w$, $\overline{F}=\phi F$, $\overline{\xi}=\phi \xi$, $\overline{U}=\phi U$ and $\overline{P}=e^{\overline{U}}$. These maps agree with the original maps on $B^4$ and are zero on $\partial B^4_2$. Another feature of these new maps is that their various Sobolev or Lorentz norms are estimated by the corresponding norms of the original maps. With the help of the equations \eqref{z5}, \eqref{z22}, \eqref{z23} and \eqref{z251} we get new maps $\overline{\omega}$, $\overline{\Omega}$ and $\overline{K}_i$, $i\in \{1,\dots,5\}$ from all these modified functions. By these considerations we see that instead of solving \eqref{z6} we can also solve the equation
\begin{equation}
\label{z253}
\begin{array}{rl}
\ds \nabla \Delta \overline{A} +\Delta \overline{A} ~ \overline{K}_2+\nabla^2\overline{A}~ \overline{K}_3- \nabla \overline{A}~ \overline{K}_4 +\overline{A}~ \overline{K}_5=(curl B) \overline{P},
\end{array}
\end{equation}
for $\overline{A}$ and $B$. First we want to solve the system
\begin{equation}
\label{z254}
\begin{array}{rl} 
\ds \Delta^2 \hat{A} &= (curl B)\cdot \nabla \overline{P} -\on{div}( \Delta \hat{A}\overline{K}_2 + \nabla^2 \hat{A}\overline{K}_3 -\nabla \hat{A} \overline{K}_4+\hat{A} \overline{K}_5+\overline{K}_5), \\[5mm]
\ds  curl (curl B)&= curl \big( ( \nabla \Delta \hat{A}) +\Delta \hat{A} \overline{K}_2 + \nabla^2 \hat{A}\overline{K}_3 +\nabla \hat{A} \overline{K}_4+\hat{A} \overline{K}_5+\overline{K}_5)\overline{P}^{-1}\big) , \\[5mm]
\ds dB&=0, \\[5mm]
\ds \hat{A}&= \frac{\partial \Delta \hat{A}}{\partial \nu}=0 \ \ \ \text{on}\ \ \ \partial B_2^4, \\[5mm]
\ds  B\cdot \nu&=0\ \ \ \text{on}\ \ \ \partial B_2^4,  \\[5mm]
\ds \int_{B_2^4} \Delta \hat{A} &=0,
\end{array}
\end{equation}
where we use the notation:
\newline
If $B=B_{kl}\ \partial_{x_k} \wedge \partial_{x_l}$ is a two-vectorfield, then 
\begin{equation}
\label{z254a}
\begin{array}{rl}
\ds dB &= \sum_{i,k,l} \frac{\partial B_{kl}}{\partial_{x_i}}\ \partial_{x_i} \wedge \partial_{x_k} \wedge \partial_{x_l},\\[5mm]
\ds B\cdot \nu&=\sum_{k,l} B_{kl}\ \nu_l\ \partial_{x_k}
\end{array}
\end{equation}
and for a (one-)vectorfield $C=C_k\ \partial_{x_k}$ we define
\begin{equation}
\label{z254b}
\begin{array}{rl}
\ds curl C &= \sum_{k,l} \Big(\frac{\partial C_k}{\partial_{x_l}}-\frac{\partial C_l}{\partial_{x_k}}\Big) \partial_{x_k} \wedge \partial_{x_l}.
\end{array}
\end{equation}
Using Lemma \ref{linfinitybounds} and the above remarks we get
\begin{equation}
\label{z255}
\begin{array}{rl} 
\ds ||\hat{A}||_{L^\infty}+||\hat{A}||_{W^{2,2}}+||\nabla \Delta \hat{A}||_{L^{\frac{4}{3},1}} \le& c||d B||_{L^{\frac{4}{3},1}}||\nabla \overline{P}^{-1}||_{L^{4,2}}\\[5mm]
\ds &+c\Delta (||\hat{A}||_{W^{2,2}}+||\hat{A}||_{L^\infty})+c\Delta.
\end{array}
\end{equation}
Furthermore from Lemma \ref{estimatesB} we get
\begin{equation}
\label{z256}
\begin{array}{rl} 
||d B||_{L^{\frac{4}{3},1}}\le& c||\nabla \Delta \hat{A}||_{L^{\frac{4}{3},1}}+c\Delta (||\hat{A}||_{W^{2,2}}+||\hat{A}||_{L^\infty})+c\Delta. 
\end{array}
\end{equation}
By a standard fixed-point argument we get the existence of solutions $\hat{A}$ and $B$. The maps $A^*=\hat{A}+\on{id}$ and $B$ then satisfy $|| \on{dist}(A^\star, SO(n))||_{L^\infty} \le c \varepsilon$ and 
\begin{equation}
\label{z257}
\begin{array}{rl} 
\nabla \Delta A^* +\Delta A^*\overline{K}_2 +\nabla^2 A^*\overline{K}_3 -\nabla A^* \overline{K}_4+A^* \overline{K}_5-(curlB) \overline{P}=~curl C,
\end{array}
\end{equation}
where $C \in W^{1,\frac{4}{3},1}(B_2^4, M(m)\otimes \wedge^2 \R^4)$. From this we see that 
\begin{equation}
\label{z258}
\begin{array}{rl}
curl(curlC \overline{P}^{-1})&=0\ \ \ \text{in} \ \ \ B_2^4, \\[5mm]
C\cdot \nu&=0\ \ \ \text{on} \ \ \ B_2^4. 
\end{array}
\end{equation}
By Lemma \ref{Czero} we get that $C$ is identically zero and therefore we get our desired solution $A$ and $B$ of equation \eqref{z6} in $B^4$.

\subsection{Proof of Theorem \ref{th-1}.}
From Theorem \ref{th-2} and Theorem \ref{th-3} we see that 
\begin{equation}
\label{z25a}
\Delta (A \Delta u)=\on{div} K,
\end{equation}
in $B_{\frac{1}{2}}$, where $K \in L^2 \cdot W^{1,2}(B_{\frac{1}{2}})\subset L^{\frac{4}{3},1}(B_{\frac{1}{2}})$. Using $L^p$-theory in Lorentz spaces we get $A\Delta u \in L_{\on{loc}}^{2,1}(B_{\frac{1}{4}})$ and therefore $\Delta u \in L_{\on{loc}}^{2,1}(B_{\frac{1}{4}})$. Let us first assume that $u$ is smooth. Then we extend $u|_{B_{\frac{1}{4}}}$ to all of $\R^4$ with compact support such that the extension $\tilde{u}$ satisfies
\begin{equation}
\label{z25b}
||\Delta \tilde{u}||_{L^{2,1}(\R^4)}+||\tilde{u}||_{W^{1,2}(\R^4)}\le c(||\Delta u||_{L^{2,1}(B_{\frac{1}{4}})}+||u||_{W^{1,2}(B_{\frac{1}{4}})})\le c.
\end{equation}
In the following we let $G(x)=c\on{log} |x|$ be the fundamental solution of $\Delta^2$ on $\R^4$. Since $|\Delta G|= O(\frac{1}{|x|^2})\in L^{2,\infty}(\R^4)$ we conclude that for all $y\in B_{\frac{1}{8}}$
\begin{equation}
\label{z25c}
\begin{array}{rl}
\ds |u(y)|=|\tilde{u}(y)| &= |\int_{\R^4} \Delta^2 G (y-z) \tilde{u}(z)dz|\\[5mm]
\ds &= |\int_{\R^4} \Delta G (y-z) \Delta \tilde{u}(z)dz|\\[5mm]
\ds &\le c||\Delta \tilde{u}||_{L^{2,1}(\R^4)}\\[5mm]
\ds &\le c(||\Delta u||_{L^{2,1}(B_{\frac{1}{4}})}+||u||_{W^{1,2}(B_{\frac{1}{4}})}).
\end{array}
\end{equation}
Using the density of smooth maps in $\{v\in W^{2,2}(B_{\frac{1}{4}}) | \ \ ||\Delta v||_{L^{2,1}(B_{\frac{1}{4}})}+||v||_{W^{1,2}(B_{\frac{1}{4}})}< \infty\}$ we see that $u \in C^0(B_{\frac{1}{8}})$. This proves Theorem \ref{th-1}. 
\subsection{Proof of Theorem \ref{th-4}.}
Smooth extrinsic biharmonic maps satisfy the equation (see for example \cite{lamm03b})
\begin{equation}
\label{z26}
\begin{array}{rl}
\ds -\Delta^2 u =& \sum_{i=n+1}^m(\Delta\langle \nabla u,(d\nu_i \circ u)\nabla u \rangle+\nabla \cdot \langle \Delta u, (d\nu_i \circ u)\nabla u \rangle \\[5mm]
\ds &+\langle \nabla \Delta u, (d\nu_i \circ u)\nabla u \rangle)\nu_i\circ u, 
\end{array}
\end{equation}
where $\{ \nu_i\}_{i=n+1}^m$ is an orthonormal frame of the normal space of $N$ near $u(x)$ for all $x\in M$. Now we rewrite this equation term by term.
\begin{equation}
\label{z27}
\begin{array}{rl}
\ds \langle \nabla \Delta u, (d\nu_i \circ u)\nabla u \rangle \nu_i\circ u &= \nabla u^k \nabla \Delta u^j( (\nu_i\circ u)^s (d_k\nu_i \circ u)^j-(\nu_i\circ u)^k (d_s\nu_i \circ u)^j)
\end{array}
\end{equation}
For the second term we have
\begin{equation}
\label{z28}
\begin{array}{rl}
\nabla \cdot \langle \Delta u, (d\nu_i \circ u)\nabla u \rangle \nu_i\circ u =& \nabla \cdot(\nabla u^k \Delta u^j( (\nu_i\circ u)^s (d_k\nu_i \circ u)^j)\\[5mm]
\ds &-\nabla u^k \Delta u^j(d_k\nu_i \circ u)^j\nabla (\nu_i\circ u)^s.
\end{array}
\end{equation}
Finally we calculate
\begin{equation}
\label{z29}
\begin{array}{rl}
\ds \Delta\langle \nabla u,(d\nu_i \circ u)\nabla u \rangle \nu_i\circ u =& \Delta (\nabla u^k \nabla u^j( (\nu_i\circ u)^s (d_k\nu_i \circ u)^j)\\[5mm]
\ds &-2 \nabla \cdot (\nabla u^k \nabla u^j (d_k\nu_i \circ u)^j \nabla (\nu_i\circ u)^s)\\[5mm]
\ds &+\nabla u^k \nabla u^j(d_k\nu_i \circ u)^j\Delta (\nu_i\circ u)^s.
\end{array}
\end{equation}
Therefore the claim follows by defining
\begin{equation}
\label{extbih}
\begin{array}{rl}
\ds W =& \nabla \Delta u^j\Big( (\nu_i\circ u)^s (d_k\nu_i \circ u)^j-(\nu_i\circ u)^k (d_s\nu_i \circ u)^j\Big)\\ [5mm]
\ds &+\nabla u^j(d_k\nu_i \circ u)^j\Delta (\nu_i\circ u)^s-\Delta u^j(d_k\nu_i \circ u)^j\nabla (\nu_i\circ u)^s,\\[5mm]
\ds w =&  \Delta u^j (\nu_i\circ u)^s (d_k\nu_i \circ u)^j -2  \nabla u^j (d_k\nu_i \circ u)^j \nabla (\nu_i\circ u)^s,\\[5mm]
V =&  \nabla u^j (\nu_i\circ u)^s (d_k\nu_i \circ u)^j,\\[5mm]
\ds \omega =& \Delta u^j\Big( (\nu_i\circ u)^s (d_k\nu_i \circ u)^j-(\nu_i\circ u)^k (d_s\nu_i \circ u)^j\Big),\\[5mm]
\ds F=&-\Delta u^j \nabla \Big( (\nu_i\circ u)^s (d_k\nu_i \circ u)^j-(\nu_i\circ u)^k (d_s\nu_i \circ u)^j\Big)\\[5mm]
\ds &+\nabla u^j(d_k\nu_i \circ u)^j\Delta (\nu_i\circ u)^s-\Delta u^j(d_k\nu_i \circ u)^j\nabla (\nu_i\circ u)^s.
\end{array}
\end{equation}
In the general case (when we can not localize in the target to find a local orthonormal frame of the normal space) we use the following equivalent equation for extrinsic biharmonic maps 
\begin{equation}
\label{extbih1}
\begin{array}{rl}
\ds \Delta^2 u= \Delta(A(u)(\nabla u ,\nabla u))+\nabla \cdot (\langle \nabla (P(u)), \Delta u\rangle )+\langle \nabla (P(u)), \nabla \Delta u\rangle,
\end{array}
\end{equation}
which was first derived by Wang \cite{wang2}, and where $A$ is the second fundamental form of the embedding $N\hookrightarrow \R^m$ and $P(y):\R^m \rightarrow T_yN$ is the orthogonal projection on the tangent space of $N$ at the point $y$. With the help of this equation we see that we can always rewrite extrinsic biharmonic maps in the desired form. 
\newline
In the case of intrinsic biharmonic maps we note that we have the relation
\begin{equation}
\label{z30}
\begin{array}{rl}
\mathcal{B}_{int}(u)=\mathcal{B}_{ext}(u)-\int_{B^4} |A(u)(\nabla u, \nabla u)|^2dx.
\end{array}
\end{equation}
This implies that the Euler-Lagrange equation of $\mathcal{B}_{int}$ differs from the Euler-Lagrange equation of $\mathcal{B}_{ext}$ only by the term coming from the variation of $\int_{B^4} |A(u)(\nabla u, \nabla u)|^2dx$. It is easy to see that this term is of the form 
\begin{equation}
\label{z30a}
\begin{array}{rl}
\nabla \cdot( g_2 \nabla u)+g_1 \nabla u,
\end{array}
\end{equation}
where $g_2 \in L^2$ and $g_1 \in L^2 \cdot W^{1,2}$ (for details see \cite{wang2}). This proves the claim.
\section{Proof of Theorem \ref{flow1}.}
\setcounter{equation}{0}
In this section we apply the conservation law obtained in section $2$ to study the existence of a unique global weak solution of the extrinsic biharmonic map flow in the energy space. 
\newline
Before proving Theorem \ref{flow1} we need a regularity result for $L^2$-perturbations of extrinsic biharmonic maps. 
\begin{lemma}\label{pertubation}
Let $f \in L^2(M,\R^m)$ and let $u\in W^{2,2}(M,N)$ be a weak solution of
\begin{align}
\Delta^2 u &=\Delta (V\cdot \nabla u)+\text{div} (w\nabla u)+W\cdot \nabla u+f\ \ \ \text{in}\ \ \ M,  \label{perequation}
\end{align}
where $V$, $w$ and $W$ are as in \eqref{extbih}. There exists $\varepsilon>0$ and $C>0$ such that if
\begin{align}
\kappa(u;R):=\int_{B_{32R}(x)}( |\nabla^2 u|^2+\frac{1}{R^2} |\nabla u|^2) < \varepsilon, \label{small}
\end{align}
for some $R>0$ and $x\in M$, then $u\in W^{4,2}(B_R(x),N)$ and we have the estimate
\begin{align}
\int_{B_R(x)}|\nabla^4 u|^2 \le \frac{c\kappa(u;R)}{R^4}+c\int_{B_{32R}(x)}|f|^2.\label{estper}
\end{align}
\end{lemma}
\begin{proof}
Let $\varepsilon$ be as in Theorem \ref{th-3}. By Theorem \ref{th-2} and Theorem \ref{th-3} we know that there exists $A \in L^\infty \cap W^{2,2}(B_{16R}(x),Gl(m))$ with $\nabla \Delta A\in L^{\frac{4}{3}}(B_{16R}(x),M(m))$, 
\newline
$||\on{dist}(A,SO(m))||_{L^\infty} \le c \varepsilon$  and $B\in W^{1,\frac{4}{3}}(B_{16R}(x),M(m)\otimes \wedge^2 \R^4)$ such that
\begin{equation}
\label{consper}
\begin{array}{rl}
\ds Af=& \Delta (A\Delta u)-\on{div} \Big( 2\nabla A \Delta u- \Delta A \nabla u -Aw \nabla u-\nabla A(V \cdot \nabla u)  \\[5mm]
\ds&+ A \nabla (V\cdot \nabla u)+ B\cdot \nabla u \Big).
\end{array} 
\end{equation}
Arguing as in the proof of Lemma \ref{linfinitybounds} we get that $\Delta u \in W^{1,\frac{4}{3}}(B_{8R}(x),\R^m)$. Next we choose a smooth cut-off function $\varphi \in C^\infty_0(B_{8R}(x))$ with compact support such that $\varphi \equiv 1$ in $B_{4R}(x)$ and $||\nabla^j \varphi||_{L^\infty} \le \frac{c}{R^j}$ for $0\le j \le 4$. From this it is easy to see that $v:=\varphi u$ solves 
\begin{equation}
\label{consper1}
\begin{array}{rl}
\ds \varphi Af+g=& \Delta (A\Delta v)-\on{div} \Big( 2\nabla A \Delta v- \Delta A \nabla v -Aw \nabla v-\nabla A(V \cdot \nabla v)  \\[5mm]
\ds &+ A \nabla (V\cdot \nabla v)+ B\cdot \nabla v \Big), 
\end{array} 
\end{equation}
where 
\begin{equation}
\label{defg}
\begin{array}{rl}
\ds g=& \Delta (A\on{div}(u\nabla \varphi))-\on{div} \Big( 2\nabla A \on{div}(u\nabla \varphi)- \Delta A u\nabla \varphi -Aw u\nabla \varphi-\nabla A(V \cdot u\nabla \varphi)  \\[5mm]
\ds &+ A \nabla (V\cdot u\nabla \varphi)+ B\cdot u\nabla \varphi \Big)\\[5mm]
\ds \in& L^{\frac{4}{3}}(B_{8R}(x)).
\end{array} 
\end{equation}
Now we define $\tilde{v}=A\Delta v$ and we claim that for $\varepsilon$ small enough and all $\frac{4}{3}\le p <2$ the operator 
\begin{equation}
\label{operator}
\begin{array}{rl}
 &\mathcal{S}:W^{1,p}(B_{8R}(x),\R^m)\rightarrow W^{1,p}(B_{8R}(x),\R^m), \\[5mm]
 &\mathcal{S}(\tilde{v})= \tilde{v}- \Delta^{-1} \on{div}\Big( 2\nabla A \Delta v-\Delta A \nabla v - Aw\nabla v-\nabla A(V\cdot \nabla v)\\[5mm]
&+A \nabla (V\cdot \nabla v)+ B\cdot \nabla v \Big)
\end{array} 
\end{equation}
is a bijection. First we consider the case $\frac{4}{3}<p<2$. We note that by the Sobolev embedding the assumption $\tilde{v}\in W^{1,p}(B_{8R}(x),\R^m)$, $\frac{4}{3}<p<2$, implies that $\nabla^2 v \in L^q(B_{8R}(x),\R^m)$ with $\frac{1}{q}+\frac{1}{4}=\frac{1}{p}$ and $\nabla v \in L^r(B_{8R}(x),\R^m)$ with $\frac{1}{r}+\frac{1}{2}=\frac{1}{p}$. Moreover we have the estimate
\begin{align}
||\nabla \Delta v||_{L^p}+||\nabla^2 v||_{L^q}+||\nabla v||_{L^r} \le c ||\tilde{v}||_{W^{1,p}}. \label{estv}
\end{align}
We estimate
\begin{equation}
\label{estper1}
\begin{array}{rl}
\ds || \on{div}\Big(& 2\nabla A \Delta v-\Delta A \nabla v - Aw\nabla v-\nabla A(V\cdot \nabla v) +A \nabla (V\cdot \nabla v)+ B\cdot \nabla v \Big)||_{W^{-1,p}} \\[5mm]
\ds \le& c||\Big( 2\nabla A \Delta v-\Delta A \nabla v - Aw\nabla v-\nabla A(V\cdot \nabla v) +A \nabla (V\cdot \nabla v)+ B\cdot \nabla v \Big)||_{L^p} \\[5mm]
\ds \le& c||\nabla A||_{L^4}(||\nabla^2 v||_{L^q}+||V||_{L^4}||\nabla v||_{L^r})  +c||B||_{L^2}||\nabla v||_{L^r} \\[5mm]
\ds &+c||\nabla^2 A||_{L^2}||\nabla v||_{L^r}+c||w||_{L^2}||\nabla v||_{L^r} +c||V||_{L^4}||\nabla^2 v||_{L^q}+c||\nabla V||_{L^2}||\nabla v||_{L^r}\\[5mm]
\ds \le& c\varepsilon ||\tilde{v}||_{W^{1,p}},
\end{array} 
\end{equation}
where we used the estimates of Theorem \ref{th-3}, \eqref{extbih} and \eqref{estv} in the last step. Since $\Delta^{-1}$ is a continuous map from $W^{-1,p}(B_{8R}(x),\R^m)$ into $W_0^{1,p}(B_{8R}(x),\R^m)$, this shows that $\mathcal{S}:W^{1,p}(B_{8R}(x),\R^m)\rightarrow W^{1,p}(B_{8R}(x),\R^m)$ is a bijection for every $\frac{4}{3}<p<2$ and $\varepsilon$ small enough.
\newline
Arguing as in the proof of Lemma \ref{linfinitybounds} it is easy to see that the same is true in the case $p=\frac{4}{3}$. This implies in particular that $\tilde{v}=A \Delta v \in W_0^{1,p}(B_{8R}(x),\R^m)$, for every $p<2$, which then implies that $u\in W^{3,p}(B_{4R}(x),\R^m)$ for every $p<2$. Next we use \eqref{perequation} iteratively to get first $u\in W^{4,p}(B_{2R}(x),\R^m)$ for every $p<2$ and then $u\in W^{4,2}(B_{R}(x),\R^m)$. 
\newline
To prove \eqref{estper} we use the interpolation inequality
\begin{equation}
\label{interpolation}
\begin{array}{rl}
\ds \int_{M}\varphi^4 (|\nabla u|^8+|\nabla^2 u|^4+|\nabla^3 u|^{\frac{8}{3}})\le c\kappa(u;R)\int_M\varphi^4 |\nabla^4 u|^2+\frac{c\kappa(u;R)}{R^4},
\end{array} 
\end{equation}
see \cite{lamm03b}, where $\varphi$ is a cut-off function as above. From this and the explicit form of the equation \eqref{extbih} the claim follows.
\end{proof}
Now we are able to prove Theorem \ref{flow1}.
\newline
\it{Proof of Theorem \ref{flow1}.} \rm
\newline
By previous results of Gastel \cite{gastel04} and Wang \cite{wang06} we know that for every $u_0\in W^{2,2}(M,N)$ there exists a unique weak solution of \eqref{flow} which is smooth away from finitely many times and for which the energy $\mathcal{B}_{ext}$ is monotonically decreasing. Let us denote the first singular time of the solution by $T=T(u_0)$ and let us denote the unique smooth solution by $v\in C^\infty(M \times (0,T),N)$. As in the work of Freire \cite{freire95} on the harmonic map flow we realize that it only remains to prove that $u=v$ in $M\times (0,T)$, where $u\in H^1([0,T],L^2(M)) \cap L^2([0,T],W^{2,2}(M))$ is a solution of the flow in the energy space. The uniqueness then implies that $u$ is smooth on $M\times (0,T)$ and we can iterate the argument.
\newline
In the following we let $U=u-v$. By Lemma \ref{pertubation} we know that $U(t) \in W^{4,2}(M,N)$ for almost every $t\in [0,T)$. Covering $M$ by balls $B_{R_t}(x_i)$ such that at most finitely many of the balls $B_{32R_t}(x_i)$ intersect and that we have 
\begin{align}
K(t,R_t):=\kappa(u(t);R_t)=\frac{\varepsilon}{2}, \label{coveri}
\end{align}
where $\varepsilon$ is as in Lemma \ref{pertubation}. A consequence of \cite{gastel04} and \cite{lamm03b} is that $T$ is characterized by the fact that $R_t \rightarrow 0$ as $t\rightarrow T$. Using the finite covering property and the estimate \eqref{estper} we get
\begin{equation}
\label{estflow}
\begin{array}{rl}
\ds \int_M |\nabla^4 u|^2 &\le \frac{cK(t,R_t)}{R_t^4}+c\int_M |\partial_t u|^2\\ [5mm]
\ds &\le c
\end{array} 
\end{equation}
for almost every $t\in (0,T)$. Therefore we can directly follow the uniqueness proof of Gastel \cite{gastel04} for the smooth situation to conclude that $U\equiv 0$ and this finishes our proof.
\qed
\appendix
\section{Appendix}
\renewcommand{\theequation}{A.\arabic{equation}}
\setcounter{equation}{0}
In the appendix we collect certain existence results and estimates for second and fourth order systems which we need in Section $2$ and $3$. Instead of proving these results for vectorfields we prove them for forms. The two-vectorfield $B=B_{kl}\ \partial_{x_k} \wedge \partial_{x_l}$ can be identified with the two-form $B=B_{kl}\ dx_k \wedge dx_l$. The differential operator $curl$ for two-vectorfields is the same as $\delta$ for two-forms and the operator $d$ for two-vectorfields is the same as the exterior derivative $d$ for two-forms. Moreover the operator $curl$ for one-vectorfields $C=C_k\ \partial_{x_k}$ is the same as $d$ for the corresponding one-form $C=C_k\ dx_k$. With the help of these identifications one easily sees that the Theorems below imply the results which were needed in the proofs of the Theorems \ref{th-1}.--\ref{th-3}.     
\begin{lemma}\label{estimatesB}
Let $f \in L^{\frac{4}{3},1}(B^4,M(m)\otimes \wedge^1 \R^4)$ and let $B \in W^{1,\frac{4}{3},1}(B^4,M(m)\otimes \wedge^2 \R^4)$ be a solution of 
\begin{equation}
\label{equationB}
\begin{array}{rl}
\ds d \delta  B &=  d f\ \ \ \text{in}\ \ \ B^4, \\[5mm]
\ds d B&=0 \ \ \ \text{in}\ \ \ B^4\ \ \ \text{and}  \\[5mm]
\ds i^\star_{\partial B^4} (\star B)&=0
\end{array} 
\end{equation}
then we have
\begin{align}
||d B||_{L^{\frac{4}{3},1}(B^4)} \le c ||f||_{L^{\frac{4}{3},1}(B^4)}. \label{estimateB}
\end{align}
\end{lemma}
\begin{proof}
Using the Hodge decomposition (see \cite{iwaniec01}, Corollary $10.5.1$) and an interpolation argument (see \cite{hunt66}) we can write
\begin{align*}
f= d\alpha +\delta \beta,
\end{align*}
where $\alpha \in W^{1,\frac{4}{3},1}(B^4,M(m))$ and $\beta \in W^{1,\frac{4}{3},1}(B^4,M(m)\otimes \wedge^2 \R^4)$ satisfy
\begin{align}
||\alpha||_{W^{1,\frac{4}{3},1}(B^4)}+||\beta||_{W^{1,\frac{4}{3},1}(B^4)} \le c ||f||_{L^{\frac{4}{3},1}(B^4)}. \label{estalphabeta}
\end{align}
Moreover we have $ i^\star_{\partial B^4} (\star \beta) =0$ and $\beta$ can be written as $\beta = d \tilde{\beta}$, where 
\newline
$\tilde{\beta}\in W^{1,\frac{4}{3},1}(B^4,M(m)\otimes \wedge^1 \R^4)$. Defining $B=\beta$ we see that $B$ solves \eqref{equationB} and the estimate \eqref{estimateB} follows from \eqref{estalphabeta}. From the fact that the homogeneous problem corresponding to \eqref{equationB} has only the trivial solution we get the desired result.  
\end{proof}
\begin{lemma}\label{Czero}
There exists $\varepsilon >0$ such that for every $P\in W^{1,4}(B^4,SO(m))$ satisfying
\begin{align}
||d P||_{L^4(B^4)}+||d P^{-1}||_{L^4(B^4)}\le \varepsilon, \label{smallnessP}
\end{align}
the only solution $C\in W^{1,\frac{4}{3}}(B^4,M(m)\otimes \wedge^2 \R^4)$ of 
\begin{equation}
\label{equC}
\begin{array}{rl}
\ds d(\delta C P^{-1})&=0,\\[5mm]
\ds i^\star_{\partial B^4}(\star C)&=0 
\end{array}
\end{equation}
is $C=0$.
\end{lemma}
\begin{proof}
First of all we claim that we can without loss of generality assume that $C=d \gamma$ for some $\gamma \in W^{1,\frac{4}{3}}(B^4,M(m)\otimes \wedge^1 \R^4)$. If this is not the case we use the Hodge decomposition for $\delta C$ (see \cite{iwaniec01}) to get
\begin{equation}
\label{equC2}
\begin{array}{rl}
\ds \delta C&= d\alpha +\delta \beta, \\[5mm]
\ds i^\star_{\partial B^4}(\star \beta) &=0\ \ \ \text{and}\\[5mm]
\ds \beta&= d \gamma.
\end{array}
\end{equation}
From this we see that $d\alpha$ is harmonic and $ i^\star_{\partial B^4}d \alpha=0$ which yields $d \alpha=0$. Therefore replacing $C$ with $\beta$ proves the claim. This implies that additional to \eqref{equC} we can assume $d C =0$.
\newline
Since $d(\delta C P^{-1})=0$ we get the existence of $f\in W^{1,\frac{4}{3}}(B^4,M(m))$ such that
\begin{align}
\delta C P^{-1} = d f. \label{equf}
\end{align}
This implies that 
\begin{equation}
\label{f}
\begin{array}{rl}
\ds \Delta f&= \star (d(\star C)\wedge dP^{-1}) \ \ \ \text{in}\ \ \ B^4,   \\[5mm]
\ds  f &=const. \ \ \ \text{on}\ \ \ \partial B^4.
\end{array}
\end{equation}
By subtracting the constant we can assume that $f$ is zero on the boundary. Using the results of Coifman, Lions, Meyer and Semmes \cite{coifman93} gives
\begin{align}
||df||_{L^{\frac{4}{3}}(B^4)}\le c ||dP^{-1}||_{L^4(B^4)} ||dC||_{L^{\frac{4}{3}}(B^4)}. \label{est1}
\end{align}
On the other hand with the remarks from above we have that
\begin{equation}
\label{equC1}
\begin{array}{rl}
\ds d \delta  C &=d(df P), \\[5mm]
\ds d C &=0,  \\[5mm]
\ds i^\star_{\partial B^4}(\star C)&=0. 
\end{array}
\end{equation}
Using Lemma \ref{estimatesB} (without the interpolation argument to get estimates in Lorentz spaces) we get
\begin{align}
||dC||_{L^{\frac{4}{3}}(B^4)}\le c ||df||_{L^{\frac{4}{3}}(B^4)}. \label{est2}
\end{align}
Combining this with \eqref{est1} we get that $d(\star C) = 0$. Using \eqref{equC1} this implies that $\star C$ is harmonic and with the help of $i^\star_{\partial B^4}(\star C)=0$ we conclude that $\star C\equiv 0$.
\end{proof}
In the next Lemma we prove a version of Wente's inequality for fourth order systems. The proof of this Lemma relies on the fact that the fundamental solution $G$ of $\Delta^2$ on $\R^4$ is a multiple of $\on{log} r$.
\begin{lemma}\label{linfinitybounds}
Let $u\in W^{2,2}( B^4,\R^m)$ be a solution of 
\begin{equation}
\label{wente}
\begin{array}{rl}
\ds-\Delta^2 u &= \star (d v \wedge d w)+\on{div} f\ \ \ \text{in}\ \ \ B^4,\\[5mm]
\ds u=\frac{\partial}{\partial \nu}\Delta u&=0\ \ \ \text{on} \ \ \partial B^4\ \ \ \text{and}\\[5mm]
\ds \int_\Omega \Delta u &=0,
\end{array}
\end{equation}
where $ v\in W^{1,p}(B^4, \R^m\otimes \wedge^2 \R^4)$, $w\in W^{1,q}(B^4,\R^m)$ ($\frac{1}{p}+\frac{1}{q}=1$) and $f\in L^{\frac{4}{3},1}(B^4 ,\R^{4m})$. Then we have
\begin{align}
||u||_{L^\infty}+||u||_{W^{2,2}}+||\Delta u||_{L^{2,1}}+||\nabla \Delta u||_{L^{\frac{4}{3},1}}\le c ||d v||_{L^p}||d w||_{L^q}+c||f||_{L^{\frac{4}{3},1}}. \label{estlinfinity}
\end{align}
\end{lemma}
\begin{proof}
First we consider a solution $u_1$ of
\begin{equation}
\label{wente1}
\begin{array}{rl}
\ds -\Delta^2 u_1 &= \star (d v \wedge d w)\ \ \ \text{in}\ \ \ B^4,\\[5mm]
\ds u_1=\frac{\partial}{\partial \nu}\Delta u_1&=0\ \ \ \text{on} \ \ \partial B^4\ \ \ \text{and}\\[5mm]
\ds \int_\Omega \Delta u_1 &=0.
\end{array}
\end{equation}
Due to the above mentioned fact that $G(x)=\on{log} |x|$ we can follow the proof of the Wente inequality for the two-dimensional case word by word (see for example \cite{brezis85} or \cite{helein02}) and get the desired estimate for $u_1$. The $L^{\frac{4}{3},1}$-estimate for $\nabla \Delta u_1$ and the $L^{2,1}$-estimate for $\Delta u_1$ follow from the work of Coifman, Lions, Meyer and Semmes \cite{coifman93} applied to $\Delta (\Delta u_1)$ and the embeddings $W^{1,1} \hookrightarrow L^{\frac{4}{3},1}$ and $W^{2,1} \hookrightarrow L^{2,1}$ in four dimensions (see \cite{poornima83}). 
\newline
Next we consider the solution $u_2$ of 
\begin{equation}
\label{wente2}
\begin{array}{rl}
\ds -\Delta^2 u_2 &=\on{div} f\ \ \ \text{in}\ \ \ B^4,\\[5mm]
\ds u_2=\frac{\partial}{\partial \nu}\Delta u_2&=0\ \ \ \text{on} \ \ \partial B^4\ \ \ \text{and}\\[5mm]
\ds \int_\Omega \Delta u_2 &=0.
\end{array}
\end{equation}
By classical $L^p$-theory and interpolation (see \cite{hunt66}) we obtain $u_2\in W^{2,2}$ and $\nabla \Delta u_2 \in L^{\frac{4}{3},1}$ with the desired estimate. To derive the $L^\infty$-bound of $u_2$ we assume first that $u_2 \in C^\infty(\overline{B^4},\R^{m})$. Then we extend $u_2$ to all of $\R^4$ with compact support and we let $\hat{u}_2$ denote this extension. Moreover, by interpolation, we can assume that $||\nabla \Delta \hat{u}_2||_{L^{\frac{4}{3},1}}\le c(||\nabla \Delta u_2||_{L^{\frac{4}{3},1}}+ ||u_2||_{W^{2,2}})$. Due to the fact that $\nabla G \in L^{4,\infty}$ we can estimate
\begin{equation}
\label{wente3}
\begin{array}{rl}
\ds |\hat{u}_2(x)|&=|\int_{\R^4} \Delta^2 G(y-x) \hat{u}_2(y)dy| \\[5mm]
\ds &=|\int_{\R^4}  \nabla G(y-x) \nabla \Delta \hat{u}_2(y)dy|\\[5mm]
\ds &\le c||\nabla \Delta \hat{u}_2||_{L^{\frac{4}{3},1}}\\[5mm]
\ds &\le c(||\nabla \Delta u_2||_{L^{\frac{4}{3},1}}+ ||u_2||_{W^{2,2}}).
\end{array}
\end{equation}
By density the $L^\infty$-bound of $u_2$ follows. Finally the $L^{2,1}$-bound for $\Delta u$ follows from the embedding $L_1^{\frac{4}{3},1}\hookrightarrow L^{2,1}$ (see \cite{oneil63}, \cite{peetre63}). Altogether this gives the desired estimates for $u$.
\end{proof}
\begin{bem}
A special situation of the above result is the case $f\in W^{1,1}(B^4, \R^m) \hookrightarrow L^{\frac{4}{3},1}(B^4,\R^m)$.
\end{bem}
Next we prove the gauge transformation result which we need for the proof of Theorem \ref{th-3}.
\begin{theorem}\label{Uhlenbeckgauge}
There exists $\varepsilon >0$ and $c>0$ such that for every $\Omega \in W^{1,2}(B^4,so(m)\otimes \wedge^1 \R^4)$ satisfying
\begin{align}
\int_{B^4}|\nabla \Omega |^2 + (\int_{B^4}|\Omega |^4)^{\frac{1}{2}} < \varepsilon, \label{smallnesgauge}
\end{align}
there exists $\xi \in W^{2,2}(B^4,so(m)\otimes \wedge^2 \R^4)$ and $U\in W^{2,2}(B^4,so(m))$ such that
\begin{equation}
\label{Uhlenbeckgauge1}
\begin{array}{rl}
\ds \Omega &= e^{-U} d e^U+e^{-U}d^\star \xi e^U,  \\[5mm]
\ds d(i_{\partial B^4}^\star \star \xi) &=0,  \\[5mm]
\ds ||U||_{W^{2,2}(B^4)}+||\xi||_{W^{2,2}(B^4)} &\le c (||\nabla \Omega||_{L^2(B^4)}+||\Omega||_{L^4(B^4)}). 
\end{array}
\end{equation}
\end{theorem}
\begin{proof}
The result follows by a compactness argument from the following Lemma \ref{gauge3}.
\end{proof}
\begin{lemma}\label{gauge3}
There exists $\varepsilon >0$ and $c>0$ such that for every $\alpha>0$ and every $\Omega \in W^{1,2+\alpha}(B^4,so(m)\otimes \wedge^1 \R^4)$ satisfying
\begin{align}
\int_{B^4}|\nabla \Omega |^2 + (\int_{B^4}|\Omega |^4)^{\frac{1}{2}} < \varepsilon, \label{smallnesgauge2}
\end{align}
there exists $\xi \in W^{2,2+\alpha}(B^4,so(m)\otimes \wedge^2 \R^4)$ and $U\in W^{2,2+\alpha}(B^4,so(m))$ such that
\begin{align}
\ds \Omega &= e^{-U} d e^U+e^{-U}d^\star \xi e^U, \label{gequ} \\[5mm]
\ds d(i_{\partial B^4}^\star \star \xi) &=0, \label{gboun} \\[5mm]
\ds ||U||_{W^{2,2}(B^4)}+||\xi||_{W^{2,2}(B^4)} &\le c (||\nabla \Omega||_{L^2(B^4)}+||\Omega||_{L^4(B^4)})\ \ \ \text{and} \label{gest1}\\[5mm]
\ds ||U||_{W^{2,2+\alpha}(B^4)}+||\xi||_{W^{2,2+\alpha}(B^4)} &\le c (||\nabla \Omega||_{L^{2+\alpha}(B^4)}+||\Omega||_{L^{4+2\alpha}(B^4)}). \label{gest2}
\end{align}
\end{lemma}
\begin{proof}
Since the proof of this result follows closely the arguments given in \cite{riviere06}, \cite{rivierestruwe06} and \cite{uhlenbeck82} we only sketch the main ideas here. For $\alpha>0$ we introduce the set
\begin{align*}
\mathcal{U}^\alpha_{\varepsilon,C}:= \{& \Omega \in W^{1,2+\alpha}(B^4,so(m)\otimes \wedge^1 \R^4)|\int_{B^4}|\nabla \Omega |^2 + (\int_{B^4}|\Omega |^4)^{\frac{1}{2}} < \varepsilon\\
&\text{and there exist} \  P   \ \text{and}  \ \xi  \ \text{satisfying}\ \ \ \eqref{gequ}-\eqref{gest2}\}.
\end{align*}
Since $\Omega \equiv 0 \in \mathcal{U}^\alpha_{\varepsilon,C}$ it remains to show that, for $\varepsilon>0$ small enough and $C$ large enough, the set $\mathcal{U}^\alpha_{\varepsilon,C}$ is both open and closed in the path-connected set\begin{align*}
\mathcal{V}^\alpha_\varepsilon:= \{ \Omega \in W^{1,2+\alpha}(B^4,so(m)\otimes \wedge^1 \R^4)|\int_{B^4}|\nabla \Omega |^2 + (\int_{B^4}|\Omega |^4)^{\frac{1}{2}} < \varepsilon \}.
\end{align*}
The proof of the closedness relies on the fact that we have the embedding $W^{2,2+\alpha}(B^4)\hookrightarrow C^{0,\beta}(B^4)$, for some $\beta>0$. This allows to pass to the limit in \eqref{gequ}. The openness follows as in \cite{riviere06}, \cite{rivierestruwe06} and \cite{uhlenbeck82}.
\end{proof}

\end{document}